\documentclass{amsart}

\def\RR{\mathbb R}
\def\CC{\mathbb C}
\newcommand{\remove}[1]{ }

\providecommand{\abs}[1]{\lvert#1\rvert}

\newcommand{\set}[1]{\left\{#1\right\}}

\newtheorem{theorem}{Theorem}[section]

\newtheorem{lemma}[theorem]{Lemma}

\theoremstyle{definition}

\theoremstyle{remark}
\newtheorem*{remark}{Remark}

\begin{document}
\title{Properties of expansions in complex bases}
\author{Vilmos Komornik}
\address{D\'epartement de math\'ematique\\
         Universit\'e de Strasbourg\\
         7 rue Ren\'e Descartes\\
         67084 Strasbourg Cedex, France}
\email{komornik@math.u-strasbg.fr}
\author{Paola Loreti}
\address{Dipartimento di Metodi e Modelli\\
         Matematici per le Scienze Applicate\\
         Sapienza Universit\`a di Roma\\
         Via A. Scarpa, 16\\
         00161 Roma, Italy}
\email{loreti@dmmm.uniroma1.it}

\begin{abstract}
Expansions in noninteger positive bases have been intensively investigated since the pioneering works of R\'enyi (1957) and Parry (1960). The discovery of surprising unique expansions in certain noninteger bases by Erd\H os, Horv\'ath and Jo\'o (1991) was followed by many studies aiming to clarify the topological and combinatorial nature of the sets of these bases. In the present work we extend some of these studies to more general, negative or complex bases. 
\end{abstract}

\maketitle


\section{Introductions}
\label{s0}

Given a noninteger real number $q>1$ in his seminal paper \cite{Ren} R\'enyi initiated the study of expansions of the form
\begin{equation*}
x=\sum_{k=1}^{\infty}\frac{c_k}{q^k}
\end{equation*}
on the digit set $\set{0,1,\ldots,[q]}$. Following a subsequent work of Parry \cite{Par60} it turned out that the properties of these expansions are radically different in many aspects from the usual expansions in integer bases. Sometimes a number may have infinitely, even continuum many different expansions, but for any given positive integer there are many real numbers having exactly that number of expansions. 

Many properties depend on whether $q$ is smaller or greater than the golden mean. It is obvious that $x=1$ has only one expansion in base $2$. It was natural to expect that $x=1$ has infinitely many different expansions in each base $1<q<2$. However, the situation is quite complex, and partly still non understood. If $q$ is smaller than the golden mean, than $x=1$ continuum many different expansions, but only countably many ones if  $q$ is equal to the golden mean. But the most surprising phenomenon was discovered by Erd\H os, Horv\'ath and Jo\'o \cite{ErdHorJoo1991} who exhibited continuum many bases $1<q<2$ in which the expansion of $x=1$ is unique. Subsequently the set $\mathcal{U}$ of such bases was characterized in \cite{ErdJooKom55} and their nature was clarified in \cite{KomLor95} and \cite{KomLor125}. For example, we have the following properties:

\begin{itemize}
\item The set $\mathcal{U}$ has zero Lebesgue measure, but it contains continuum many elements.

\item The set $\mathcal{U}$ is closed from above and it has a smallest element $q\approx 1.787$.

\item The closure of $\mathcal{U}$ is a Cantor set (i.e., a closed set having no isolated or interior points) of zero Lebesgue measure, and  $\overline{\mathcal{U}}\setminus \mathcal{U}$ is a countable dense set in $\overline{\mathcal{U}}$.
\end{itemize}

If $q$ is less than the golden mean, then not only $x=1$ but all real numbers satisfying $0<x<1/(q-1)$ have continuum many different expansion. This raises the question whether there exist expansion with special properties, for instance universal expansions which contain all possible finite variations of the digits zero and one. It was proven in \cite{ErdKom87} that if $q$ sufficiently close to one, then every real number satisfying $0<x<1/(q-1)$ has universal expansions.

Expansions in complex bases were also investigated by Dar\'oczy and I. K\'atai \cite{DarKat1986} and then in \cite{KomLor136}. 
The purpose of this paper is to extend some of these results to negative or even complex bases. 

In the following by a \emph{sequence} we always mean a sequence of  zeroes and ones. Given a real or complex number $q$ of modulus $\abs{q}>1$, by an \emph{expansion} of a real or complex number $x$ in base $q$ we mean a sequence $(c_k)$ satisfying
\begin{equation*}
\sum_{k=1}^{\infty}\frac{c_k}{q^k}=x.
\end{equation*}
We denote by $J_q$ the set of numbers having at least one expansion in base $q$.

\section{Expansions in positive bases}
\label{s1}

Fix a positive real number $q>1$. Then $J_q\subset\left[ 0,\frac{1}{q-1}\right] $ with equality if and only if $q\le 2$; for $q>2$ the set $J_q$ is not connected. The endpoints of $J_q$ always have the unique expansions $0^{\infty}$ and $1^{\infty}$. 

We recall from \cite{ErdHorJoo1991}, \cite{ErdJooKom55}, \cite{ErdKom87}, \cite{KomLor95} and \cite{KomLor125} the following results:

\begin{theorem}\label{t11}\mbox{}

(a) If $q=\frac{1+\sqrt{5}}{2}$, then the set of expansions of $x=1$ is countably infinite.
\smallskip

(b) If $q<\frac{1+\sqrt{5}}{2}$, then each interior point of $J_q$ has continuum many expansions.

(c) If, moreover, $q<\sqrt[4]{2}$ and $q^2$ is different from the second Pisot number, then each interior point of $J_q$ has a universal expansion, i.e., an expansion containing all possible finite blocks of zeroes and ones.
\end{theorem}

\section{Expansions in negative bases}
\label{s2}

In this section we fix a negative real number $q=-p$ with $p>1$. The identity
\begin{equation*}
\sum_{k=1}^{\infty}\frac{c_k}{q^k}
=\frac{-p}{p^2-1}+\sum_{k=1}^{\infty}\frac{1-c_{2k-1}}{p^{2k-1}}+\frac{c_{2k}}{p^{2k}}
\end{equation*}
establishes a bijection between the expansions in bases $q$ and $p$: \emph{$(c_k)$ is an expansion of $x$ in base $q$ if and only if $(1-c_1,c_2,1-c_3,c_4,\ldots)$ is an expansion of $x+\frac{p}{p^2-1}$ in base $p$.} In particular, we have
\begin{equation*}
J_q=\frac{-p}{p^2-1}+J_p.
\end{equation*}
Using the results of the preceding section it follows that
\begin{equation*}
J_q\subset\left[ \frac{q}{q^2-1},\frac{1}{q^2-1}\right] =\left[ \frac{-p}{p^2-1},\frac{1}{p^2-1}\right]
\end{equation*}
with equality if and only if $p\le 2$; for $p>2$ the set $J_q$ is not connected. The endpoints of $J_q$  always have the unique expansions $(10)^{\infty}$ and $(01)^{\infty}$.

Let us observe that the sequences $(c_k)$ and $(1-c_1,c_2,1-c_3,c_4,\ldots)$ are universal at the same time. Indeed, this is the case $m'=1$ of the following lemma:

\begin{lemma}\label{l21}
Fix a positive integer $m'$ and for each sequence $d=(d_i)$  of zeroes and ones let us define another sequence $d'=(d_i')=T(d)$ of zeroes and ones by the formula
\begin{equation}\label{transform}
d'_{m'i+j}:=
\begin{cases}
1-d_{m'i+j}&\text{for $m'=0,2,\ldots$ and $j=1,\ldots, m'$,}\\
d_{m'i+j}&\text{for $m'=1,3,\ldots$ and $j=1,\ldots, m'$.}
\end{cases}
\end{equation}
Then $d$ and $d'$ are universal at the same time.
\end{lemma}

\begin{proof}
First we show that if $d$ is universal, then $d'$ contains any fixed finite block $B'$. We may assume that the length of $B'$ is a multiple of $m'$. It will be convenient to define the transform $T$ by the same formula \eqref{transform} for finite sequences as well.

By the universality of $d$ it contains the finite block
\begin{equation*}
B:=T(B')0^{m-1}T(B')0^{m-2}T(B')\ldots T(B')00T(B')0T(B')
\end{equation*}
and, depending on the position of $B$ in $d$ one of these $m$ blocks $T(B')$ will be transformed into $T(T(B'))=B'$ in $d'=T(d)$, so that $d'$ contains the block $B'$ indeed. (If the index of the first element of the block $B$ in $d$ is $m'i+j$ for some $1\le j\le m'$, then the 
$j$th subblock $T(B')$ will be transformed into a block $B'$.) Since $T(T(d))=d$ for all sequences, by symmetry the universality of $d'$ implies that of $d$, too.
\end{proof}

We deduce from Theorem \ref{t11} and Lemma \ref{l21} the following results:

\begin{theorem}\label{t22}\mbox{}

(a) If $q=-\frac{1+\sqrt{5}}{2}$, then the set of expansions of $x=0$ is countably infinite.
\smallskip

(b) If $-\frac{1+\sqrt{5}}{2}<q<-1$, then each interior point of $J_q$ has continuum many expansions.
\smallskip

(c) If, moreover, $-\sqrt[4]{2}<q<-1$ and $q^2$ is different from the second Pisot number, then each interior point of $J_q$ has a universal expansion, i.e., an expansion containing all possible finite blocks of zeroes and ones.
\end{theorem}

\section{Expansions in purely imaginary bases}
\label{s3}

If $q=ip$ for some $p>1$, then by distinguishing the even and odd powers of $q$ we obtain that $J_q=qJ_{q^2}+J_{q^2}$. Since $q^2=-p^2$, identifying $\CC$ with $\RR^2$ we deduce from the results of the preceding section that
\begin{equation*}
J_q\subset\left[ \frac{-p^2}{p^4-1},\frac{1}{p^4-1}\right]\times  \left[ \frac{-p^3}{p^4-1},\frac{p}{p^4-1}\right]
\end{equation*}
with equality if and only if $\abs{q}\le \sqrt{2}$. The four vertices of this closed rectangle always have the unique expansions $(1000)^{\infty}$, $(0100)^{\infty}$, $(0010)^{\infty}$ and $(0001)^{\infty}$.

Furthermore, we deduce from Theorem \ref{t22} the following result:

\begin{theorem}\label{t31}\mbox{}

(a) If $\abs{q}=\sqrt{\frac{1+\sqrt{5}}{2}}$, then the set of expansions of $x=0$ is countably infinite.
\smallskip

(b) If $\abs{q}<\sqrt{\frac{1+\sqrt{5}}{2}}$, then each point of $J_q$ except the four vertices has continuum many expansions.
\end{theorem}

\begin{remark}
Analogous results hold if $q=-ip$ for some $p>1$.
\end{remark}

The existence of universal expansions follows from the case $m=4$ of the more general Theorems \ref{t42} and \ref{t44} in Section \ref{s4} below.

\section{Complex bases}
\label{s4}

Set $\CC_R:=\set{z\in\CC\ :\ \abs{z}\le R}$ for brevity. 

\begin{theorem}\label{t41}
Fix a complex number $\omega$ of modulus one such that $\omega^4\ne 1$ and a positive 
real number $R$.  If $p>1$ is sufficiently 
close to $1$, then every complex number 
$z\in\CC_R$ has continuum many expansions in base $q:=p\omega$.
\end{theorem}

\begin{proof}
Applying Proposition 2.1 from \cite{KomLor136} with $\omega^2$ and $2R$ instead of $\omega$ and $R$ we obtain that every complex number 
$w\in\CC_{2R}$ has at least one expansion in base $q^2$. Now fix a sufficiently large positive integer $n$ such that
\begin{equation*}
\frac{1}{p^{2n+1}}+\frac{1}{p^{2n+3}}+\cdots <R.
\end{equation*}
Now given $z\in\CC_R$ and an arbitrary subset $A\subset\set{2n+1,2n+3,\ldots}$, applying this result with
\begin{equation*}
w:=z-\sum_{k\in A}\frac{1}{q^k}\in \CC_{2R}
\end{equation*}
we obtain an expansion
\begin{equation*}
z-\sum_{k\in A}\frac{1}{q^k}=\frac{c_2}{q^2}+\frac{c_4}{q^4}+\cdots .
\end{equation*}
This yields a different expansion of $z$ in base $q$ for every subset $A$.
\end{proof}

Now fix a positive integer $m\ge 3$, a complex number satisfying $\omega^m=1$ and consider $q=p\omega$ with $1<p\le 2^{1/m}$. It follows from the identity
\begin{equation}\label{due}
\sum_{i=1}^{\infty}\frac{d_i}{q^i}=\sum_{j=1}^m(p\omega)^{m-j}\sum_{i=0}^{\infty}\frac{d_{mi+j}}{p^{m(i+1)}}
\end{equation}
and from the inequalities $1<p^m\le 2$ that $J_q$ coincides with the set of numbers
\begin{equation*}
\sum_{j=1}^m(p\omega)^{m-j}\alpha_j
\end{equation*}
with  
\begin{equation*}
0\le \alpha_j\le \frac{1}{p^m-1}\quad\text{for all}\quad j.
\end{equation*}

We will prove the following theorem:

\begin{theorem}\label{t42}\mbox{}
If $1<p<2^{1/4m}$ and $p^{2m}$ is different from the second Pisot number, then each interior point of $(p^m-1)J_q$ has a universal expansion, i.e., an expansion containing all possible finite blocks of zeroes and ones.
\end{theorem}

We divide the proof into three steps.
\medskip

\emph{First step.} 
Every interior point $z$ of $(p^m-1)J_q$ may be written in the form
\begin{equation*}
z=\sum_{j=1}^m(p\omega)^{m-j}\alpha_j
\end{equation*}
with suitable real numbers satisfying $0< \alpha_j<1$ for all $j$.
Indeed, since $tz\in (p^m-1)J_q$ for a sufficiently small $t>1$, dividing by $t$ a representation of $tz$ we may assume that $0\le \alpha_j<1$ for all $j$. Adding a small positive multiple of the equality
\begin{equation*}
0=\sum_{j=1}^m \omega^{m-j}
\end{equation*}
we may also assume that $0<\alpha_j<1$ for all $j$.

Using the identity \eqref{due} it suffices to construct a universal sequence $(d_i)$ satisfying the equalities
\begin{equation}\label{41a}
\alpha_j=\sum_{i=0}^{\infty}\frac{d_{mi+j}}{p^{m(i+1)}},\quad j=1,\ldots, m.\qedhere
\end{equation} 

\emph{Second step.} 
The construction will be based on the following generalization of Lemma 4.1 in \cite{ErdKom87}.

\begin{lemma}\label{l43}
Let $m$ and $p$ be as in the statement of the theorem. 
Given $0<\alpha_j\le 1$,  $j=1,\ldots, m$ and a finite sequence $c_1\ldots c_{mN}$ of integers $\in\{0,1\}$ there exists another finite sequence $d_1\ldots d_{mN+mn}$ of integers $\in\{0,1\}$ ending with $c_1\ldots c_{mN}$:
\begin{equation*}
d_{mn+i}=c_i,\quad i=1,\ldots, mN,
\end{equation*}
and satisfying the inequalities
\begin{equation}\label{41}
0<\alpha_j-\sum_{i=0}^{n+N-1}\frac{d_{mi+j}}{p^{m(i+1)}}<\frac{1}{p^{m(n+N)}},\quad j=1,\ldots, m.
\end{equation}
\end{lemma}

\begin{proof}
Let $0=y_0<y_1<\cdots$ be the sequence of numbers of the form $P(p^m)$ where $P$ runs over the polynomials with coefficients belonging to the set $\{0,1\}$. We know from \cite{ErdKom87} that $y_k\to\infty$ and $y_{k+1}-y_k\to 0$. Set
\begin{equation*}
A_j=\sum_{i=0}^{N-1}\frac{c_{mi+j}}{p^{m(i+1)}},\quad j=1,\ldots, m.
\end{equation*}

If $n$ is sufficiently large, then
\begin{equation*}
p^{mn}\alpha_{j}>A_{j},\quad j=1,\ldots, m,
\end{equation*}
so that there exist indices $k_1,\ldots,k_m$ satisfying
\begin{equation*}
y_{k_j}<p^{mn}\alpha_j-A_j\le y_{k_j+1},\quad j=1,\ldots, m.
\end{equation*}
Since $y_{k_j}<p^{mn}$ for all $j$, we have
\begin{equation*}
y_{k_j}=a_{j,0}+a_{j,m}p^m+\cdots+a_{j,(n-1)m}p^{(n-1)m},\quad j=1,\ldots, m,
\end{equation*}
with suitable coefficients $a_{j,km}\in\{0,1\}$. 

Since $n\to\infty$ implies $k_j\to\infty$ for each $j$, choosing a sufficiently large $n$ we also have
\begin{equation*}
0<y_{k_j+1}-y_{k_j}<p^{-mN},\quad j=1,\ldots, m.
\end{equation*}
It follows that
\begin{equation*}
0<\alpha_j-p^{-mn}y_{k_j}-p^{-mn}A_j<p^{-m(n+N)},\quad j=1,\ldots, m.\qedhere
\end{equation*}
\end{proof}
\medskip

\emph{Third step.} 
We fix a sequence $B_1, B_2,\ldots$ of finite blocks, containing all possible finite variations of the digits $0$ and $1$. Applying the preceding lemma with $c_1\ldots c_{mN}:=B_1$ we obtain a finite sequence $d_1,\ldots, d_{mi_1}$ (with $i_1:=n+N$) ending with the block $B_1$ and satisfying the inequalities
\begin{equation*}
0<\alpha_j-\sum_{i=0}^{i_1-1}\frac{d_{mi+j}}{p^{m(i+1)}}<\frac{1}{p^{mi_1}},\quad j=1,\ldots, m.
\end{equation*}

Repeating the construction with the numbers
\begin{equation*}
\alpha_{j,2}:=p^{mi_1}\left( \alpha_j-\sum_{i=0}^{i_1-1}\frac{d_{mi+j}}{p^{m(i+1)}}\right),\quad j=1,\ldots, m
\end{equation*}
instead of $\alpha_j$ and with the block $B_2$ instead of $B_1$, we may extend the sequence $d_1,\ldots, d_{mi_1}$ to a  sequence $d_1,\ldots, d_{mi_2}$  ending with the block $B_2$ and satisfying the inequalities
\begin{equation*}
0<\alpha_{j,2}-\sum_{i=i_1}^{i_2-1}\frac{d_{mi+j}}{p^{m(i+1-i_1)}}<\frac{1}{p^{m(i_2-i_1)}},\quad j=1,\ldots, m.
\end{equation*}
This can be rewritten in the following form:
\begin{equation*}
0<\alpha_j-\sum_{i=0}^{i_2-1}\frac{d_{mi+j}}{p^{m(i+1)}}<\frac{1}{p^{mi_2}},\quad j=1,\ldots, m.
\end{equation*}

Proceeding by induction we construct in this way a sequence $(d_i)$ and an infinite sequence of indices $i_1<i_2<\cdots$ such $d_1,\ldots, d_{mi_k}$ ends with the block $B_k$, and 
\begin{equation*}
0<\alpha_j-\sum_{i=0}^{i_k-1}\frac{d_{mi+j}}{p^{m(i+1)}}<\frac{1}{p^{mi_k}},\quad j=1,\ldots, m
\end{equation*}
for each $k=1,2,\ldots .$, so that  the equalities \eqref{41a} are satisfied.

If $m$ is even, then we may weaken the assumptions on $p$ in order to obtain universal expansions.

\begin{theorem}\label{t44}\mbox{}
Assume that $m=2m'$ is even, $1<p<2^{1/2m}$ and $p^m$ is different from the second Pisot number. Then every number of the form
\begin{equation}\label{41uj}
z=\sum_{j=1}^{m'}(p\omega)^{m'-j}\alpha_j
\end{equation}
with real numbers
\begin{equation*}
-\frac{p^{m'}}{p^{2m'}-1}<\alpha_j<1-\frac{p^{m'}}{p^{2m'}-1}=\frac{p^{2m'}-p^{m'}-1}{p^{2m'}-1}
\end{equation*}
has a universal expansion.
\end{theorem}

\begin{proof}
It follows from the identity
\begin{align*}
\sum_{i=1}^{\infty}\frac{d_i}{q^i}
&=\sum_{i=1}^{\infty}\frac{d_i}{(p\omega)^i}\\
&=\sum_{j=1}^{m'}(p\omega)^{m'-j}\sum_{i=0}^{\infty}\frac{d_{m'i+j}}{(p\omega)^{m'(i+1)}}\\
&=\sum_{j=1}^{m'}(p\omega)^{m'-j}\sum_{i=0}^{\infty}(-1)^{i+1}\frac{d_{m'i+j}}{(p^{m'})^{i+1}}
\end{align*}
and from the results of Section \ref{s2} (recall that $1<p^{2m'}\le 2$) that a complex number $z$ belongs to $J_q$ if and only if $z$ has the form \eqref{41uj}
with suitable real numbers
\begin{equation}\label{42uj}
\alpha_j\in \left[ \frac{-p^{m'}}{p^{2m'}-1},\frac{1}{p^{2m'}-1}\right].
\end{equation}
Rewriting the above identity in the form
\begin{multline}\label{sei}
\left( \sum_{i=1}^{\infty}\frac{d_i}{q^i}\right) +\frac{p^{m'}}{p^{2m'}-1}\sum_{j=1}^{m'}(p\omega)^{m'-j} \\
=\sum_{j=1}^{m'}(p\omega)^{m'-j}\left( \frac{1-d_j}{p^{m'}}+\frac{d_{j+m'}}{(p^{m'})^2}+\frac{1-d_{j+2m'}}{(p^{m'})^3}+\frac{d_{j+3m'}}{(p^{m'})^4}-\cdots\right) 
\end{multline}
we also see that $z\in J_q$ if and only if 
\begin{equation*}
z':=z+\frac{p^{m'}}{p^{2m'}-1}\sum_{j=1}^{m'}(p\omega)^{m'-j}=\sum_{j=1}^{m'}(p\omega)^{m'-j}\alpha_j
\end{equation*}
with suitable real numbers
\begin{equation*}
\alpha_j\in \left[ 0,\frac{1}{p^{m'}-1}\right].
\end{equation*}

Furthermore, the identity \eqref{sei} establishes a bijection between the expansions 
\begin{equation}\label{expansionz}
z=\sum_{i=1}^{\infty}\frac{d_i}{q^i}=\sum_{j=1}^{m'}(p\omega)^{m'-j}\sum_{i=0}^{\infty}(-1)^{i+1}\frac{d_{m'i+j}}{(p^{m'})^{i+1}}
\end{equation}
of an element $z\in J_q$ and the expansion
\begin{equation}\label{expansionz'}
z'=\sum_{i=1}^{\infty}\frac{d_i'}{q^i}=\sum_{j=1}^{m'}(p\omega)^{m'-j}\sum_{i=0}^{\infty}\frac{d_{m'i+j}'}{(p^{m'})^{i+1}}
\end{equation}
of its translate
\begin{equation}\label{translate}
z':=z+\frac{p^{m'}}{p^{2m'}-1}\sum_{j=1}^{m'}(p\omega)^{m'-j},
\end{equation}
where the sequence $(d_i')$ of zeroes and ones is obtained from $(d_i)$ by the transformation \eqref{transform} of Lemma \ref{l21}.

If $z$ satisfies the hypotheses of the theorem, then by what we have shown in the first step, its translate $z'$ defined by \eqref{translate} may be written in the form
\begin{equation*}
z'=\sum_{j=1}^{m'}(p\omega)^{m'-j}\alpha_j
\end{equation*}
with  real numbers $0<\alpha_j<1$. Repeating the second and third part of the proof of Theorem \ref{t42} we conclude that $z'$ has a universal expansion \eqref{expansionz'}. By Lemma \ref{l21} its transform $d:=T(d')$ provides a universal expansion \eqref{expansionz}. 
\end{proof}

\end{document}